\documentstyle[12pt,amscd,times]{amsart}
\input amssymb.sty
\newtheorem{theorem}{Theorem}[section]
\newtheorem{lemma}{Lemma}[section]

\newtheorem{corollary}{Corollary}[section]
\newtheorem{remark}{Remark}[section]

\theoremstyle{definition}
\newtheorem{definition}{Definition}[section]

\newcommand{\frakt}{\mathfrak}

\def\ad{\mathop{\operatorname {ad}}\nolimits}

\def\Ext{\mathop{\operatorname{{\cal E}xt}}\nolimits}
\def\Fred{\mathop{\operatorname {Fred}}\nolimits}

\def\Hom{\mathop{\operatorname {Hom}}\nolimits}

\def\Ind{\mathop{\operatorname{Ind}}\nolimits}

\def\Ker{\mathop{\operatorname{Ker}}\nolimits}
\def\Lie{\mathop{\operatorname {Lie}}\nolimits}

\def\Mat{\mathop{\operatorname {Mat}}\nolimits}

\def\Nil-Rad{\mathop{\operatorname {Nil-Rad}}\nolimits}

\def\rank{\mathop{\operatorname {rank}}\nolimits}

\def\semi-linear{\mathop{\operatorname {semi-linear}}\nolimits}
\def\sgn{\mathop{\operatorname {sgn}}\nolimits}

\def\SU{\mathop{\operatorname {SU}}\nolimits}
\def\SO{\mathop{\operatorname {SO}}\nolimits}
\def\sgn{\mathop{\operatorname{sgn}}\nolimits}

\def\ker{\mathop{\operatorname{Ker}}\nolimits}
\def\HP{\mathop{\operatorname {HP}}\nolimits}
\def\HE{\mathop{\operatorname {HE}}\nolimits}
\def\Tot{\mathop{\operatorname {Tot}}\nolimits}

\begin{document}
\title[Group C*-Algebras \& Compact Quantum
Groups]{On the Structure of Group C*-Algebras and Compact Quantum
Groups\footnote{This is the text of the ICM98 Talk - Section 7 ``Lie
Groups and Lie Algebras", 
Saturday, August 22, 1998}} \author{Do Ngoc Diep}
\maketitle
\section{Motivation }

Consider a locally compact group $G$ and consider some
appropriate group algebras. The group algebras ${\Bbb C}[G]$ for $G$ as an 
abstract group is not enough to define the structure of $G$. We must find a
more effective group algebra. 

For any locally compact group $G$ 
there is a natural left-(right-)invariant Haar measure $dg$. 

The space
$L^2(G):= L^2(G,dg)$ of the square-integrable functions plays an important 
role in harmonic analysis. If the group is of type I, $L^2(G)$ admits a 
spectral decomposition with respect to
 the left and right regular representations into a sum of the direct sum 
(the so called discrete series) and/or the direct integral (the continuous 
series) of irreducible unitary representations. 

The space $L^1(G)=L^1(G,dg)$ 
of the  functions with integrable module plays a crucial role. 
With the well-defined convolution product, 
$$\varphi,\psi \in L^1(G) \mapsto \varphi * \psi \in L^1(G);$$
$$(\varphi * \psi)(x) := \int_G \varphi(y)\psi(y^{-1}x)dy
$$ it becomes a Banach algebra.
There is also a well-defined Fourier-Gel'fand transformation on $L^1(G)$,
$$\varphi \in L^1(G,g) \mapsto \hat{\varphi},$$
$$\hat{\varphi}(\pi) := \pi(\varphi) = \int_G \pi(x)\varphi(x)dx.$$

There is a one-to-one correspondence between the (irreducible) unitary 
representations of $G$ and the  non-degenerate (irreducible) *-representations 
of the involutive Banach algebra $L^1(G)$. The general theorems of the spectral theory of the 
representations of $G$ are then proved with the help of an appropriate 
translation in to the corresponding theory of $L^1(G)$, for which one can use more 
tools from functional analysis and topology. One can also define on $L^1(G)$ 
an involution $\varphi \mapsto \varphi^*$,
$$\varphi^*(g) := \overline{\varphi(g^{-1})}.$$ 

However the norm of the involutive Banach algebra $L^1(G)$ is not regular, 
i.e. in general
$$\Vert a^*a\Vert_{L^1(G)} \ne \Vert a\Vert^2_{L^1(G)}\;.$$
It is therefore more useful to consider the corresponding regular norm 
$\Vert.\Vert_{C^*(G)}$,
$$\Vert \varphi  \Vert_{C^*(G)} := \sup_{\pi \in \hat{G}} \Vert \pi(\varphi)
\Vert$$ and its completion $C^*(G)$. The spectral theory of
unitary representations
of $G$ is equivalent to the spectral theory of non-degenerate *-representations
 of the C*-algebra $C^*(G)$. 
The general theorems of harmonic analysis say that
the structure of $G$ can be completely definite by the structure of $C^*(G)$.

One poses therefore the problem of description of 
the structure of the C*-algebras of locally compact groups. 
\begin{enumerate}
\item[1)] How to realize the irreducible unitary representations of the 
locally compact group $G$.
\item[2)] How to describe the images of the Fourier-Gel'fand transformation
and in particular, of the inclusion of $C^*(G)$ into some ``continuous'' 
product of the algebras ${\cal L}({\cal H}_\pi)$, $\pi\in \hat{G}$ of bounded
operators in the separable Hilbert space ${\cal H}_\pi$ of representation 
$\pi$.
\end{enumerate}

\subsection{Compact Groups}

Consider for the moment a compact group $G$. 
For compact group, all representations are unitarizable, 
i.e. are equivalent to some unitary ones. It is well-known also that :

\begin{enumerate}
\item[a)]
The family of irreducible representations is not more than countable.
\item[b)]
Each irreducible representation is finite dimensional, say 
$n_i, i= 1,\dots, \infty$ and there is some good realizations of these 
representations, say in tensor spaces, or last time, in cohomologies.
\item[c)]
The Fourier-Gel'fand transformation gives us an isomorphism
$$C^*(G) \cong {\prod_{i=1}^\infty}' \Mat_{n_i}({\Bbb C}),$$
where the prime in the product means the subset of "continuous vanishing at 
infinity" elements. \end{enumerate}

This means that in compact group case the group C*-algebra  plays the 
same role as the group algebra of finite groups.

\subsection{Locally Compact Group} 

{\it The main problem is how to describe the group algebra in general, and 
in particular the C*-algebra $C^*(G)$. }

We are trying to decompose our C*-algebras into some towers of ideals 
and step-by-step define the associated extensions by KK-functors or 
their generalizations. The resulting invariants form just our index. 
This idea was proposed in and developed in 
for a large class of type I C*-algebras. 

We propose a general 
construction and some reduction procedure of the K-theory invariant 
$\Ind\, C^*(G)$ of group C*-algebras. Using the orbit method, we reduces $Index~C^*(G)$ to a family of 
Connes' 
foliation C*-algebras  indices $Index~C^*(V_{2n_i},{\cal F}_{2n_i})$, by a family of KK-theory invariants. 
Using some generalization of the Kasparov type condition (treated by G.G. 
Kasparov in the nilpotent Lie group case, we reduces 
every \newline $Index C^*(V_{2n_i},{\cal F}_{2n_i})$ to a family of 
KK-theory invariants of the same type valuated in KK(X,Y) type groups. 
The last ones are in some sense computable by using the cup-cap product 
realizing the Fredholm operator indices.

\subsection{ BDF K-functor ${\cal E}xt$.}

Denote by $C(X)$ the C*-algebra of continuous complex-valued functions
over a fixed metrizable compact $X$, ${\cal H}$ a fixed separable Hilbert
space over complex numbers, ${\cal L}({\cal H})$ and ${\cal K}({\cal H})$
the C*-algebras of bounded and respectively, compact linear operators in
${\cal H}$. An extension of C*-algebras means a short exact sequence of
C*-algebras and *-homomorphisms of special type
$$0 \longrightarrow {\cal K}({\cal H}) \longrightarrow {\cal E} 
\longrightarrow C(X)\longrightarrow 0.$$ Two extensions are by definition
equivalent iff there exists an isomorphism
$ \psi : {\cal E_1} \longrightarrow {\cal E_2}$  and
its restriction
$\psi\vert_{{\cal K}({\cal H}_1)} : {\cal K}({\cal H}_1) 
\longrightarrow {\cal K}({\cal H}_2)$ such that the following diagram
is commutative
$$ 
\begin{array}{ccccccccc}
0 & \longrightarrow & {\cal K}({\cal H}_1 ) & \longrightarrow & {\cal E}_1 & \longrightarrow & C(X) & \longrightarrow & 0\\
 &    & \Big\downarrow\vcenter{%
     \rlap{$\psi\vert_.$}} & & \Big\downarrow\vcenter{\rlap{$\psi$}}
  &  & \Big\Vert &  &  \\
0 & \longrightarrow & {\cal K}({\cal H}_2) & \longrightarrow  & {\cal E}_2 
& \longrightarrow & C(X) &\longrightarrow & 0
\end{array}
$$

\subsection{The group of all affine transformations of the real
straight line}

Let us denote by $G$ the group of all affine transformations of the real
straight line. 

\begin{theorem} Every irreducible unitary representation of group $G$ is
unitarilly equivalent to one of the following mutually nonequivalent
representations: \begin{enumerate} \item[a)] the representation $S$,
realized in the space $L^2({\Bbb R}^*, \frac{dx}{ \vert x \vert})$, where
${\Bbb R}^* := {\Bbb R} \setminus (0)$, and acting in according with the
formula $$(S_gf)(x) = e^{\sqrt{-1}bx} f(ax),\text{ where } g =
\begin{pmatrix} \alpha & b \\ 0 & 1\end{pmatrix}.$$ \item[b)] the
representation $U^\varepsilon_\lambda$, realized in ${\Bbb C}^1$ and given
by the formula $$U_\lambda^\varepsilon (g) = \vert \alpha
\vert^{\sqrt{-1}\lambda}.(\sgn\alpha)^\varepsilon, \text{ where }\lambda
\in {\Bbb R}; \varepsilon = 0,1.$$ \end{enumerate} \end{theorem} 

This list of all the irreducible unitary representations gives the
corresponding list of all the irreducible non-degenerate unitary
*-representations of the group C*-algebra $C^*(G)$. It
was proved that

\begin{theorem} The group C*-algebra with formally adjoined unity
$C^*(G)^\sim$ can be included in a short exact sequence of C*-algebras and
*-homomorphisms $$ 0 \longrightarrow {\cal K} \longrightarrow C^*(G)^\sim
\longrightarrow C({\Bbb S}^1 \vee {\Bbb S}^1) \longrightarrow 0,$$ i.e.
the C*-algebra $C^*(G)^\sim$, following the BDF theory, is defined by an
element, called the index and denoted by $Index\,C^*(G)^\sim$, of the
groups ${\cal E}xt({\Bbb S}^1 \vee {\Bbb S}^1) \cong {\Bbb Z} \oplus {\Bbb
Z}$. \end{theorem}

\begin{theorem}
$$Index\, C^*(G) = (1,1) \in {\cal E}xt({\Bbb S}^1 \vee {\Bbb S}^1) 
\cong {\Bbb Z} \oplus {\Bbb Z}.$$
\end{theorem}

\section{Multidimensional Orbit Methods}

Let us consider now a connected and simply connected Lie group $G$ 
with Lie algebra ${\frakt g} := \Lie(G)$. Denote by ${\frakt g}_{\Bbb C}
$ the complexification of ${\frakt g}$. The complex conjugation in 
the Lie algebra will be also denoted by an over-line sign. 
Consider the dual space ${\frakt g}^*$ to the Lie algebra ${\frakt g}$. 
The group $G$ acts on itself by the inner automorphisms
$$A(g):= g.(.).g^{-1} : G \longrightarrow G,$$ for each $g\in G$, 
conserving the identity element $e$ as some fixed point. It follows 
therefore that the associated adjoint action $A(g)_* $ maps 
${\frakt g} = T_eG$ into itself and the co-adjoint action 
$K(g) := A(g^{-1})^* $ 
maps the dual space 
${\frakt g}^*$ 
into itself. The orbit space 
${\cal O}(G):= {\frakt g}^*/G$ 
is in general a bad topological space, namely non- Hausdorff. Consider 
one orbit 
$\Omega\in{\cal O}(G)$
and an element 
$F \in {\frakt g}^*$ in it. The stabilizer is denote by 
$G_F$, 
its connected component by $(G_F)_0$ and its Lie algebra by 
${\frakt g}_F := \Lie(G_F)$. It is well-known that
$$
\begin{array}{ccc}
G_F & \hookrightarrow & G\\
    &                 & \Big\downarrow\\
    &                  & \Omega_F
    \end{array}
$$
is a principal bundle with the structural group $G_F$. 
Let us fix some {\it connection in this principal bundle, } 
\index{connection on 
principal bundle} i.e. some {\it trivialization } \index{trivialization}
of this bundle.  
We want to construct representations in some cohomology spaces 
with coefficients in the sheaf of sections of some vector bundle 
associated with this principal bundle. 
It is well know  that every vector bundle 
is an induced one with respect to some representation of the structural group
in the typical fiber.
It is natural to fix some unitary representation 
$\tilde{\sigma}$ 
of $G_F$ such that its kernel contains $(G_F)_0$, the character
$\chi_F$ of the connected component of stabilizer
$$\chi_F(\exp{X}) := \exp{(2\pi\sqrt{-1}\langle F,X\rangle )}$$
and therefore the differential 
$D(\tilde{\sigma}\chi_F)=\tilde{\rho}$ 
is some representation  of the Lie algebra ${\frakt g}_F$. 
We suppose that the representation $D(\tilde{\rho}\chi_F)$ was extended to 
the complexification $({\frakt g}_F)_{\Bbb C}$. 
The whole space of all sections seems to be so large for the construction 
of irreducible unitary representations. 
One consider the invariant subspaces with the help of some so 
called polarizations.

\begin{definition}
We say that a triple 
$({\frakt p}, \rho,\sigma_0)$ is some 
$(\tilde{\sigma},F)$-{\it polarization, } 
iff :
\begin{enumerate}
\item[a)]
${\frakt p}$ is some subalgebra of the complex Lie algebra 
$({\frakt g})_{\Bbb C}$, containing ${\frakt g}_F$.
\item[b)]
The subalgebra ${\frakt p}$ is invariant under the action 
of all the operators of type $Ad_{{\frakt g}_{\Bbb C}}x, x\in G_F$.
\item[c)]
The vector space ${\frakt p} + \overline{\frakt p}$ 
is the complexification of some real subalgebra 
${\frakt m} = ({\frakt p} + \overline{\frakt p}) \cap {\frakt g}$. 
\item[d)]
All the subgroups $M_0$, $H_0$, $M$, $H$ are closed. where by definition 
$M_0$ (resp., $H_0$) is the connected subgroup of 
$G$ with the Lie algebra ${\frakt m}$ 
(resp., ${\frakt h} := {\frakt p} \cap {\frakt g}$) and $M:= G_F.M_0$,
$H:= G_F.H_0$.
\item[e)]
$\sigma_0$ is an irreducible representation of $H_0$ in some Hilbert space 
$V$ such that : 
1. the restriction $\sigma|_{G_F \cap H_0}$ 
is some multiple of the restriction 
$\chi_F.\tilde{\sigma}|_{G_F \cap H_0}$, 
where by definition $\chi_F(\exp X) := \exp{(2\pi\sqrt{-1}\langle F,X\rangle )}$;
2. under the action of $G_F$ on the dual $\hat{H}_0$, 
the point $\sigma_0$ is fixed.
\item[f)]
$\rho$ is some representation of the complex Lie algebra 
${\frakt p}$ in $V$ , 
which satisfies the E. Nelson conditions for $H_0$ and 
$\rho|_{\frakt h} = D\sigma_0$.
\end{enumerate}
\end{definition}

\begin{theorem}
Let us keep all the introduced above notation of $\Omega_F$, $\tilde{\sigma}$,
$G_F$,etc. and 
let us denote $\chi_F$ the character of the group $G_F$ such that 
$D\chi_F = 2\pi\sqrt{-1}F|_{{\frakt g}_F}$. Then :
\begin{enumerate}
\item[1)]
On the K-orbit $\Omega_F$ there exists a structure of some mixed manifold 
of type $(k,l,m)$, where $$k = \dim G - \dim M,$$
$$l = {1 \over 2}(\dim M - \dim H),$$
$$m = \dim H - \dim G_F.$$
\item[2)]
There exists some irreducible unitary representation 
$\sigma$ of the group $H$ such that its restriction $\sigma|_{G_F}$ 
is some multiple of the representation 
$\chi_F.\tilde{\sigma}$ and $\rho|_{\frakt h} = D\sigma$.
\item[3)]
On the $G$-fiber bundle ${\cal E}_{\sigma|_{G_F}} = G \times_{G_F} V$ 
associated with the representation $\sigma|_{G_F}$, there exists a structure
of a partially invariant and partially holomorphic Hilbert vector 
$G$-bundle ${\cal E}_{\sigma,\rho}$ 
such that the natural representation of $G$ on the space of (partially 
invariant and partially holomorphic) sections is equivalent to the 
representation by right translations of $G$ in the space 
$C^\infty (G; {\frakt p}, \rho, F, \sigma_0)$ of $V$-valued 
$C^\infty$-functions on $G$ satisfying the equations 
$$f(hx) = \sigma(h)f(x), \forall h \in H, \forall x\in G,$$
$$L_Xf + \rho(X)f = 0, \forall X\in \overline{\frakt p},$$
where $L_X$ denotes the Lie derivative along the vector field 
$\xi_X$ on $G$, corresponding to $X$.
\end{enumerate}
\end{theorem}

Consider a symplectic manifold $(M,\omega)$, 

The vector space $C^\infty (M,\omega)$, 
with respect to the Poisson brackets  
$$f_1, f_2 \in C^\infty \mapsto \{ f_1,f_2\} \in C^\infty(M,\omega)$$
become an infinite dimensional Lie algebra. 

\begin{definition} 
A {\it procedure of quantization }  \index{procedure of quantization}
is a correspondence associating to each {\it classical quantity } 
\index{classical quantity} $f\in C^\infty(M)$ a {\it quantum quantity } 
\index{quantum quantity}
$Q(f) \in {\cal L}({\cal H})$, i.e. a continuous, perhaps unbounded,
normal operator, which is auto-adjoint if $f$ is a real-valued function,
in some Hilbert space ${\cal H}$, such that
$$Q(\{f_1,f_2\}) = {i \over \hbar}[Q(f_1),Q(f_2)],$$
$$Q(1) = Id_{\cal H},$$
where $\hbar := h/2\pi$ is the normalized Planck constant, 
and $h$ is the unnormalized Planck's constant.
\end{definition}

Let us denote by ${\cal E}$ a fiber bundle into Hilbert spaces, 
$\Gamma$ a fixed connection conserving the Hilbert structure 
on the fibers; in other words,
If $\gamma$ is a curve connecting two points $x$ and $x'$, the parallel 
transport along the way $\gamma$ provides an scalar preserving isomorphism 
from the fiber ${\cal E}_x$ onto the fiber ${\cal E}_{x'}$. In this case we can
define the corresponding covariant derivative $\nabla_\xi$, $\xi \in Vect(M):=
Der\, C^\infty(M)$
in the space of smooth sections. One considers the invariant Hilbert space 
$L^2({\cal E}_{\rho,\sigma})$, which
is the completion of the space $\Gamma({\cal E}_{\rho,\sigma})$
of square-integrable partially invariant and partially holomorphic sections.

Suppose from now on that $M$ is a homogeneous $G$-space.
Choose a trivialization $\Gamma$ of the principal bundle 
$G_x \rightarrowtail G \twoheadrightarrow M$, 
where $G_x$ is the stabilizer of the fixed point $x$ on $M$.
Let us denote by $L_\xi$ the Lie derivation corresponding to the vector field
$\xi\in Vect(M)$. Let us denote by $\beta\in \Omega^1(M)$ the form of affine
connection on ${\cal E}$, corresponding to the
connection $\Gamma$ on the principal bundle. It is more comfortable to
consider the normalized connection form
$\alpha(\xi) = {\hbar \over \sqrt{-1}}\beta(\xi)$,
the values of which are anti-auto-adjoint operators on fibers. One has therefore
$$\nabla_\xi = L_\xi + {\sqrt{-1}\over \hbar}\alpha(\xi)$$

For each function $f \in C^\infty(M)$ one denotes $\xi_f$ the corresponding 
Hamiltonian vector field, i.e.
$$i(\xi_f)\omega + df = 0.$$

\begin{definition}
We define the geometrically quantized operator $Q(f)$ as
$$Q(f) := f + {\hbar\over \sqrt{-1}}\nabla_{\xi_f} = f + {\hbar\over\sqrt{-1}}
L_{\xi_f} +\alpha(\xi_f).$$
\end{definition}

\begin{theorem}

The following three conditions are equivalent.
\begin{enumerate}
\item[1)] 
$$\xi\alpha(\eta)-\eta\alpha(\xi)-\alpha([\xi,\eta])+{\sqrt{-1}\over\hbar}
[\alpha(\xi),\alpha(\eta)] = -\omega(\xi,\eta).Id; \forall \xi,\eta.$$
\item[2)]
The curvature of the affine connection $\nabla$ is equal to $-{\sqrt{-1}\over
\hbar}\omega(\xi,\eta).Id$, i.e.
$$[\nabla_\xi , \nabla_\eta] - \nabla_{[\xi,\eta]} = -{\sqrt{-1}\over\hbar}
\omega(\xi,\eta).Id; \forall \xi, \eta.$$
\item[3)]
The correspondence $f\mapsto Q(f)$ is a quantization procedure.
\end{enumerate}
\end{theorem}

Suppose that the Lie group   $G$ act  on $M$ by the symplectomorphisms. 
Then each element  $X$ of the Lie algebra ${\frakt g}$ corresponds to 
one-parameter
subgroup $\exp{(tX)}$ in $G$, which acts on $M$. Let us denote by $\xi_X$
the corresponding strictly Hamiltonian vector field. Let us denote also $L_X$
the Lie derivation along this vector field. We have 
$$[L_X,L_Y] = L_{[X,Y]},$$ and $$L_Xf = \{f_X,f\}.$$

Suppose that $f_X$ depends linearly on $X$. One has then a 2-cocycle of
the action
$$c(X,Y) := \{f_X,f_Y\} - f_{[X,Y]}.$$

\begin{definition}
We say that the action of $G$ on $M$ is {\it flat } \index{flat action}
iff this 2-cocycle is trivial.
\end{definition}

In this case we obtain from the quantization procedure a representation 
$\wedge$ of the Lie algebra ${\frakt g}$
by the anti-auto-adjoint operators 
$$X \mapsto {\sqrt{-1}\over\hbar}Q(f_X)$$
and also a representation of ${\frakt g}$ by the functions
$$X \mapsto f_X.$$ If the E. Nelson conditions are satisfied, we have a unitary
representation of the universal covering of the group $G$.

\begin{theorem}
The Lie derivative of the partially invariant and holomorphically induced
representation $\Ind(G;{\frakt p},F,\rho,\sigma_0)$ of a connected Lie 
group $G$
is just the representation obtained from the procedure of multidimensional
geometric quantization, corresponding to a fixed connection $\nabla$ of 
the partially invariant partially holomorphic induced unitarized bundle
$\overline{\cal E}_{\sigma,\rho}$, i.e.
$$\Lie_X(\Ind(G;{\frakt p},F,\rho,\sigma_0)) = {\sqrt{-1}\over\hbar}Q(f_X).$$
\end{theorem}

\section{KK-theory Invariant $Index C^*(G)$}

Let us denote by G a connected and simply connected Lie group,
${\frakt g} = \Lie (G)$ its Lie algebra, ${\frakt g}^* = \Hom_{\Bbb
R}({\frakt g}, {\Bbb R})$ the dual vector space,${\cal O} = {\cal O}%
(G)$ the space of all the co-adjoint orbits of G in ${\frakt g}^*$. This space is a disjoint union of subspaces of co-adjoint orbits of  fixed dimension, i.e. 
$${\cal O}%
= \amalg_{0 \leq 2n \leq \dim G}{\cal O}%
_{2n} ,$$
$${\cal O}%
_{2n} := \{ \Omega \in {\cal O}%
; \dim \Omega = 2n \} .$$
We define  
$$V_{2n} := \cup_{\dim \Omega = 2n}\Omega .$$
Then it is easy to see that $V_{2n}$ is the set of points of a fixed rank of the Poisson structure bilinear function
$$\{X,Y\}(F) = \langle F,[X,Y]\rangle  ,$$
 suppose it is a foliation, at least for $V_{2n}$, with $2n = max$ . 

First, we shall show that the foliation $V_{2n}$
can be obtained by the associated action of ${\Bbb R}^{2n}$ on $V_{2n}$ via 
2n times repeated action of ${\Bbb R}$ . 

Indeed, fixing any basis $X_1,X_2,\dots ,X_{2n}$ of the tangent space
${\frakt g} /{\frakt g}_F$ of $\Omega$ at the point $F \in \Omega$ ,
we can define an action ${\Bbb R}^{2n} \curvearrowright V_{2n}$ as
$$({\Bbb R} \curvearrowright ({\Bbb R} \curvearrowright (\dots {\Bbb
R} \curvearrowright V_{2n})))$$ by $$(t_1,t_2,\dots,t_{2n})
\longmapsto \exp(t_1X_1)\dots\exp(t_{2n}X_{2n}) F .$$ Thus we have the
Hamiltonian vector fields $$\xi_k := {d \over dt} |_{t=0}\exp(t_kX_k)F
, k = 1,2,\dots,2n$$ and the linear span $$F_{2n} =
\{\xi_1,\xi_2,\dots,\xi_{2n}\}$$ provides a tangent distribution.

\begin{theorem} 
$(V_{2n},F_{2n})$ is a measurable foliation. 
\end{theorem}

\begin{corollary}
 The Connes C*-algebra $C^*(V_{2n},F_{2n}) , o \leq {2n} \leq
\dim G$ are well defined.
\end{corollary}

\subsection{Reduction of $Index C^*(G)$ to $Index C^*(V_{2n},F_{2n})$ }

  Now we assume that the orbit method gives us a  complete list of irreducible representations
of $G$ , $$\pi_{\Omega_F,\sigma} = Ind(G,\Omega_F,\sigma ,{\frakt p}), \sigma \in {\cal X}%
_G(F) ,$$
the finite set of Duflo's data.

    Suppose that
$${\cal O} = \cup_{i=1}^k{\cal O}%
_{2n_i}$$ 
is the decomposition of the orbit space on a stratification of orbits of dimensions $2n_i$, where $n_1 > n_2>  \dots > n_k > 0$

We include $C^*(V_{2n_1},F_{2n_1})$ into $C^*(G)$. It is well known that the Connes C*-algebra of foliation can be included in the algebra of pseudo-differential operators of degree 0 as an ideal. This algebra of pseudo-differential operators of degree 0

is included in C*(G).

We define $$J_1 = {\bigcap_{\Omega_F \in {\cal O}%
(G) \setminus {\cal O}%
_{2n_1}}} \Ker  \pi_{\Omega_F ,\sigma},$$
and $$ {A_1} = {C^*(G)/J_1}.$$
Then $$C^*(G)/C^*(V_{2n_1},F_{2n_1}) \cong A_1$$
and we have 

$$\vbox{\halign{ #\quad &#\quad &#\quad &#\hfill\cr
$ 0 \rightarrow $ & $J_1 \rightarrow$ & $C^*(G) \rightarrow$ & $A_1 \rightarrow 0$ \cr
          \hfill      & $\enskip\downarrow$ & $\quad \downarrow Id$ & $\quad\downarrow$   \cr
$0 \rightarrow$ & $C^*(V_{2n_1},F_{2n_1}) \rightarrow$ & $C^*(G) \rightarrow$ & $C^*(G)/C^*(V_{2n_1},F_{2n_1}) \rightarrow 0$ \cr}}$$

Hence $J_1 \simeq C^*(V_{2n_1},F_{2n_1})$ and we have 
$$O \rightarrow C^*(V_{2n_1},F_{2n_1}) \rightarrow C^*(G) \rightarrow A_1 \rightarrow 0 .$$   
Repeating the procedure in replacing $$C^*(G),C^*(V_{2n_1},F_{2n_1}),A_1,J_1$$ by $$A_1,C^*(V_{2n_1},F_{2n_1}),A_2,J_2,$$ we have 
$$0 \rightarrow C^*(V_{2n_2},F_{2n_2}) \rightarrow A_1 \rightarrow A_2 \rightarrow 0$$
etc .... 

So we obtain the following result.

\begin{theorem}
 The group C*-algebra C*(G)  can be included in a finite sequence of extensions
$$ 0 \rightarrow C^*(V_{2n_1},F_{2n_1}) \rightarrow C^*(G) \rightarrow 
A_1 \rightarrow 0\leqno{(\gamma_1):}$$
$$ \quad 0 \rightarrow C^*(V_{2n_2},F_{2n_2}) \rightarrow A_1 
\rightarrow A_2 \rightarrow 0,\leqno{(\gamma_2):}$$
$$\ldots\dots\dots \dots \dots$$
$$ 0 \rightarrow C^*(V_{2n_k},F_{2n_k}) \rightarrow A_{k-1} \rightarrow A_k \rightarrow 0,\leqno{(\gamma_k):}$$
where $\widehat{A_k} \simeq Char(G)$
\end{theorem}

\begin{corollary}
$Index C^*(G)$ is reduced to the system $Index C^*(V_{2n_i},$ $ 
F_{2n_i}), i = 1,2 ,\dots , k$  by  the  invariants 
$$[\gamma_i] \in KK(A_i,C^*(V_{2n_i},F_{2n_i})) , i= 1,2,\dots ,k.$$
\end{corollary}

\begin{remark}
Ideally, all these invariants $[\gamma_i]$ could be computed step-by-step from $[\gamma_k]$ to $[\gamma_1]$.
\end{remark}

\subsection{Reduction of $Index C^*(V_{2n_i},F_{2n_i})$ to the computable
extension indices valuated in topological KK-groups of pairs of spaces}

Let us consider $C^*(V_{2n_i},F_{2n_i})$ for a fixed i. We introduce the 
following assumptions which were considered by Kasparov 
in nilpotent cases:

$(A_1)$ \quad 
There exists $k \in {\Bbb Z} , 0 < k \leq 2n_i$ such that the foliation 
$$V_{gen} := V_{2n_i} \setminus (\Lie \Gamma)^\perp$$
has its C*- algebra 
$$C^*(V_{gen},F|_{V_{gen}}) \cong C({\cal O}%
_{gen}^\sim) \otimes {\cal K}%
(H) ,$$
where $$\Gamma := {\Bbb R}^k \hookrightarrow {\Bbb R}^{2n_i} 
\hookrightarrow G ,$$  
$$\Lie \Gamma = {\Bbb R}^k \hookrightarrow {\frakt g}/{\frakt g}_{F_i} ,
(\Lie \Gamma)^\perp \subset {\frakt g}^* \cap V_{2n_i}.$$

 It is easy to see that  if the condition $(A_1)$ holds,
$C^*(V_{2n_i},F_{2n_i})$ is an extension of 
$C^*(V_{2n_i} \setminus V_{gen},F_{2n_i}|_.)$ by $C({\cal O}%
_{gen}^\sim) \otimes {\cal K}%
(H)$, where ${\cal O}_{gen}^\sim = \{ \pi_{\Omega_F,\sigma}; 
\Omega_F \in {\cal O}%
_{gen},\sigma \in {\cal X}%
_G(F)\}$, described by the multidimensional orbit method from the
previous section.
If $k = 2n_i , ({\Bbb R}^{2n_i})^\perp = \{ O \}$, 
$V_{2n_i} = V_{gen}$, we have $$C^*(V_{2n_i},F_{2n_i}) \simeq C({\cal O}%
_{2n_i}^\sim \otimes {\cal K}%
(H)) .$$
If $k = k_1 < 2n_i$ , then ${\Bbb R}^{2n_i-k_1}$ acts on $V_{2n_i} 
\setminus V_{gen}$ and we suppose that  a similar assumption $(A_2)$ holds

$(A_2)$ There exists $k_2,0 < k_2 \leq 2n_i-k_1$ such that 
$$(V_{2n_i} \setminus V_{gen})_{gen} := (V_{2n_i} \setminus V_{gen}) \setminus 
({\Bbb R}^{k_2})^\perp$$
has its C*-algebra 
$$C^*((V_{2n_i} \setminus V_{gen})_{gen},F_{2n_i}|.) \simeq C(({\cal O}%
_{2n_i} \setminus {\cal O}%
_{gen})_{gen})^\sim \otimes {\cal K}
(H) .$$
As above, if $k_2 = 2n_i - k_1$ , $C^*(V_{2n_i} \setminus V_{gen},F_{2n_i}|.) 
\simeq C(({\cal O}%
_{2n_i} \setminus {\cal O}%
_{gen})_{gen}^\sim) \otimes {\cal K}%
(H)$. In other case we repeat the procedure and go to assumption $(A_3)$, etc....

The procedure must be finished after a finite number of steps, say in m -th step,$$C^*((\dots(V_{2n_i} \setminus V_{gen}) \setminus (V_{2n_i} \setminus V_{gen})_{gen} \setminus \dots ,F_{2n_i}|_.) \simeq C((\dots ({\cal O}%
_{2n_i} \setminus {\cal O}%
_{gen}) \setminus \dots )) \otimes {\cal K}%
(H) .$$
Thus we have the following result.

\begin{theorem}
If all the arising assumptions $(A_1),(A_2),\dots$ hold, the C*-algebra $C^*(V_{2n_i},F_{2n_i})$ can be included in a finite sequence of extensions
$$ 0 \rightarrow C({\cal O}%
_{gen}^\sim) \otimes {\cal K}%
(H) \rightarrow C^*(V_{2n_i},F_{2n_i}) \rightarrow C^*(V_{2n_i} \setminus V_{gen},F_{2n_i}|_.) \rightarrow 0$$
$$ 0 \rightarrow C(({\cal O}%
_{2n_i} \setminus {\cal O}%
{gen})_{gen}^\sim) \otimes {\cal K}%
(H) \rightarrow C^*(V_{2n_i} \setminus V_{gen},F_{2n_i}) \rightarrow C^*(\dots)
\rightarrow 0$$
$$\ldots \ldots \ldots \ldots \ldots$$
$$ 0 \rightarrow C((\dots({\cal O}%
_{2n_i} \setminus {\cal O}%
_{gen}) \setminus ({\cal O}%
_{2n_i} \setminus {\cal O}%
_{gen}))_{gen}\dots{^\sim}) \otimes {\cal K}%
(H) \rightarrow $$ $$ \rightarrow C^*(\dots) \rightarrow C^*(\dots) \otimes 
{\cal K}%
(H) \rightarrow 0.$$
\end{theorem}

\subsection{General remarks concerning computation of Index C*(G)}

We see that the general computation procedure of Index C*(G) is reduced to the case of short exact sequences of type 
$$ 0 \rightarrow C(Y) \otimes {\cal K}%
(H) \rightarrow {\cal E}%
\rightarrow C(X) \otimes {\cal K}%
(H) \rightarrow 0, \leqno{(\gamma)}$$
and the index is
$$ [\gamma] = Index {\cal E}%
\in KK(X,Y) .$$
The group $KK_i(X,Y)$ can be mapped onto 
$$\oplus_{j \in {\Bbb Z}/(2)}\Hom_{\Bbb Z}(K^{i+j}(X),K^{i+j+1}(Y))$$
with kernel 
$$\oplus_{j \in {\Bbb Z}/(2)}\Ext_{\Bbb Z}^1(K^{i+j}(X),K^{i+j+1}(Y))$$
by the well known cap-product. So $[\gamma] = (\delta_0,\delta_1)$
$$\delta_0 \in \Hom_{\Bbb Z}(K^0(X),K^1(Y)) = \Ext_0(X) \wedge K^1(Y)$$
$$\delta_1 \in \Hom_{\Bbb Z}(K^1(X),K^0(Y)) = \Ext_1(X) \wedge K^0(Y).$$
Suppose $e_1,e_2,\dots,e_n \in \pi^1(X)$ to be generators and $\phi_1,\phi_2,
\dots,\phi_n \in {\cal E}$%
\quad the corresponding Fredholm operators, $T_1,T_2,\dots,T_n$ the 
Fredholm operators, representing the generators of $K^1(Y) = Index 
[Y,\Fred ]$ . We have therefore $$[\delta_0] = \sum_jc_{ij}\enskip Index T_j ,
$$ where
$$\delta_0 = (c_{ij}) \in \Mat_{\rank K^0(X) \times rank K^1(Y)}({\Bbb Z}) .$$
In the same way $\delta_1$ can be computed.

\section{Entire Periodic Cyclic Homology}
{\bf (joint with Nguyen Van Thu')}

Given a collection $\{I_\alpha\}_{\alpha\in\Gamma}$
of ideals in $A$, equipped with a so called $\ad_A$-invariant trace
$$\tau_\alpha : I_\alpha
\to {\mathbf C},$$ satisfying the properties: \begin{enumerate}
\item $\tau_\alpha$ is a {\it continuous linear} functional, normalized as
$\Vert\tau_\alpha\Vert = 1$,
\item $\tau_\alpha$ is {\it positive} in the sense that $$\tau_\alpha(a^*a)
\geq 0, \forall
\alpha \in \Gamma,$$ where the map $a \mapsto a^*$ is the involution
defining the
involutive Banach algebra structure, i.e. an anti-hermitian endomorphism
such that $a^{**}
= a$ \item $\tau_\alpha$ is {\it strictly positive} in the sense that
$\tau_\alpha(a^*a) = 0$ iff
$a=0$, for every $\alpha\in \Gamma$. \item $\tau_\alpha$ is $\ad_A$-{\it
invariant} in the
sense that $$\tau_\alpha(xa) = \tau_\alpha(ax), \forall x\in A, 
a\in I_\alpha,$$
\end{enumerate}
then we have for every $\alpha \in \Gamma$ a scalar product $$\langle
a,b\rangle_\alpha
:=\tau_\alpha(a^*b)$$ and also an inverse system $\{I_\alpha,
\tau_\alpha\}_{\alpha \in
\Gamma}$. Let $\bar{I}_\alpha$ be the completion of $I_\alpha$ with respect to the
scalar product defined
above and $\widetilde{\bar{I}_\alpha}$ denote $\bar{I}_\alpha$ with formally
adjoined unity
element. Define $C^n(\widetilde{\bar{I}_\alpha})$ as the set of
$n+1$-linear maps $\varphi :
(\widetilde{\bar{I}_\alpha})^{\otimes(n+1)} \to {\mathbf C}$. There exists a Hilbert
structure on  $(\widetilde{\bar{I}_\alpha})^{\otimes(n+1)}$ and we can identify
$C_n((\widetilde{\bar{I}_\alpha})) := \Hom(C^n(\bar{\tilde{I}}_\alpha), {\mathbf
C})$
with $C^n(\bar{\tilde{I}}_\alpha)$ via an anti-isomorphism.

For $I_\alpha \subseteq I_\beta$, we have a
well-defined map $$D^\beta_\alpha :
C^n(\widetilde{\bar{I}_\alpha}) \to C^n(\widetilde{\bar{I}_\beta}),$$ which
makes $\{
C^n(\widetilde{\bar{I}}_\alpha)\}$ into a direct system. Write $Q = \varinjlim
C^n(\widetilde{\bar{I}_\alpha})$ and observe that it admits a Hilbert
space 
structure. Let $C_n(A) := \Hom(\varinjlim
C^n(\widetilde{\bar{I}_\alpha}), {\mathbf C}) = \Hom(\varinjlim
C_n(\widetilde{\bar{I}_\alpha}), {\mathbf C})$ which is anti-isomorphic to
$\varinjlim_\alpha C^n(\bar{\tilde{I}}_\alpha)$. So we have finally
$$C_n(A) = \varinjlim C_n(\bar{\tilde{I}}_\alpha).$$

Let $$b, b' : C^n(\widetilde{\bar{I}}_\alpha) \to
C^{n+1}(\widetilde{\bar{I}}_\alpha),$$
$$N : C^n(\widetilde{\bar{I}}_\alpha) \to
C^n(\widetilde{\bar{I}}_\alpha),$$ $$\lambda :
C^n(\widetilde{\bar{I}}_\alpha) \to C^n(\widetilde{\bar{I}}_\alpha),$$
$$S : C^{n+1}(\widetilde{\bar{I}}_\alpha)\to
C^n(\widetilde{\bar{I}}_\alpha)$$ be defined as in A. Connes. Denote by $b^*, (b')^*, N^*,
\lambda^*, S^*$ the corresponding adjoint operators. Note also that for
each
$I_\alpha$ we have the same formulae for adjoint operators for homology as
Connes obtained for cohomology.

We now have a bi-complex
$$
\leqno{{\mathcal C}(A):}  
\CD
 @.    \vdots  @.    \vdots @.     \vdots @.      \\
@.          @V(-b')^* VV        @Vb^* VV       @V(-b')^* VV        @.   \\
\dots @<1-\lambda^* << \varinjlim_\alpha C_1(\bar{\tilde{I}}_\alpha) @<N^*<<
\varinjlim_\alpha C_1(\bar{\tilde{I}}_\alpha)
@<1-\lambda^*<< \varinjlim_\alpha C_1(\bar{\tilde{I}}_\alpha) @<N^*<< \dots\\
@.          @V(-b')^* VV        @Vb^* VV       @V(-b')^* VV        @.   \\
\dots @<1-\lambda^*<< \varinjlim_\alpha C_0(\bar{\tilde{I}}_\alpha) @<N^*<<
\varinjlim_\alpha 
C_0(\bar{\tilde{I}}_\alpha) @<1-\lambda^*<<\varinjlim_\alpha
C_0(\bar{\tilde{I}}_\alpha) @<N^*<< \dots\\
\endCD$$
with $d_v := b^*$ in the even columns and $d_v:= (-b')^*$ in the odd columns, $d_h
:= 1-\lambda^*$ from odd to even columns and $d_h := N^*$ from even to odd columns, 
where * means the corresponding adjoint operator. Now we have
$$\Tot({\mathcal C}(A))^{even} = \Tot({\mathcal C}(A))^{odd} :=
\oplus_{n\geq 0}
C_n(A),$$ which is periodic with period two. Hence, we have
$$\begin{array}{ccc}
&\partial & \\
\oplus_{n\geq 0} C_n(A)&\begin{array}{c} \longleftarrow\\ \longrightarrow
\end{array} &\oplus_{n\geq 0} C_n(A)\\
&\partial & \\
\end{array}$$
where $\partial = d_v + d_h$ is the total differential.

\begin{definition} Let $\HP_*(A)$ be the homology of the total complex
$(\Tot{\mathcal C}(A))$. It is called the {\sl periodic cyclic homology} 
of $A$.
\end{definition}

Note that this $\HP_*(A)$ is, in general different from the 
$\HP_*(A)$ of Cuntz-Quillen, because we used the direct limit of periodic
cyclic homology of ideals. But in special cases, when the whole
algebra $A$ is one of these ideals
 with $\ad_A$-invariant trace, (e.g. the
commutative algebras
of complex-valued functions on compact spaces), we return to the
Cuntz-Quillen $\HP_*$, which we shall use later. 

\begin{definition}
An even (or odd) chain $(f_n)_{n\geq 0}$ in ${\mathcal C}(A)$ is called
{\it entire} if the
radius of convergence of the power series $\sum_n \frac{n!}{[\frac{n}{2}]!}
\Vert
f_n\Vert z^n$, $z\in {\mathbf C}$ is infinite.
\end{definition}

Let $C^e(A)$ be the sub-complex of $C(A)$ consisting of entire chains. Then
we have a
periodic complex.

\begin{theorem} Let
$$\Tot(C^e(A))^{even} = \Tot(C_e(A))^{odd} := \bigoplus_{n\geq 0} C^e_n(A),$$
where $C^e_n(A)$ is the entire $n$-chain. Then we have a complex of
entire chains with the total differential $\partial := d_v + d_h$
$$\begin{array}{ccc} & \partial
& \\
\bigoplus_{n\geq 0} C^e_n(A) & \begin{array}{c}\longleftarrow \\
\longrightarrow\end{array} &
\bigoplus_{n\geq 0}
C^e_n(A)\\
& \partial & \end{array}$$
\end{theorem}
\begin{definition}
The homology of this complex is called also the {\it entire homology} and
denoted by $\HE_*(A)$. 
\end{definition}
Note that this entire homology is defined through
the inductive limits of ideals with ad-invariance trace. 

The main properties of this theory, namely
\begin{itemize}
\item Homotopy invariance,
\item Morita invariance and
\item Excision,
\end{itemize}
were proved and hence that $\HE_*$ is a generalized homology theory.

\begin{lemma}
If the Banach algebra $A$ can be presented as a direct limit
$\varinjlim_\alpha I_\alpha$ of a system
of ideals $I_\alpha$ with $\ad$-invariant trace $\tau_\alpha$, then
$$K_*(A) = \varinjlim_\alpha K_*(I_\alpha),$$
$$\HE_*(A) = \varinjlim_\alpha \HE_*(I_\alpha).$$
\end{lemma}

The following result from K-theory is well-known:
\begin{theorem} The entire homology of non-commutative de Rham currents
admits the following stability property
$$K_*({\mathcal K}({\mathbf H})) \cong K_*({\mathbf C}),$$
$$K_*(A \otimes {\mathcal K}({\mathbf H})) \cong K_*(A),$$
where ${\mathbf H}$ is a separable Hilbert space and $A$ is an arbitrary Banach
space.
\end{theorem}

The similar result is true for entire homology $\HE_*$ :
\begin{theorem} The entire homology of non-commutative de Rham currents
admits the following stability property
$$\HE_*({\mathcal K}({\mathbf H})) \cong \HE_*({\mathbf C}),$$
$$\HE_*(A \otimes {\mathcal K}({\mathbf H})) \cong \HE_*(A),$$
where ${\mathbf H}$ is a separable Hilbert space and $A$ is an arbitrary Banach
space.
\end{theorem}

\section{Chern Characters of Compact Lie Groups C*-Algebras}
{\bf (joint with A. O. Kuku and N. Q. Tho)}

Let $A$ be an involutive Banach algebra. In this section, we construct a
non-commutative
character
$$ch_{C^*} : K_*(A) \to \HE_*(A)$$ and later show that when $A= C^*(G)$,
this Chern character reduces up to isomorphism to classical Chern
character.

Let $A$ be an involutive Banach algebra with unity. 
\begin{theorem}
There exists a Chern character $$ch_{C^*} : K_*(A) \to \HE_*(A).$$ 
\end{theorem}

Our next result computes the Chern character in 2.1 for $A=C^*(G)$ by
reducing it to the
classical case.

\begin{theorem}
Let ${\mathbf T}$ be a fixed maximal torus of $G$ with Weyl group $W:=
N_G({\mathbf T})/{\mathbf T}$.
Then the Chern character $$ch_{C^*} : K_*({\mathbf C}^*(G)) \to
\HE_*(C^*(G))$$ is an
isomorphism, which can be identified with the classical Chern character
$$ch: K^W_*({\mathbf C}({\mathbf T})) \to \HE^W_*({\mathbf C}({\mathbf
T}))$$ that
is also an isomorphism.
\end{theorem}

\begin{remark}
For some classical groups, e. g. $\SU(n+1)$, $\SO(2n+1)$, $\SU(2n)$,
$Sp(n)$ etc. the
groups $$K^*(G) \cong K_*^W({\mathbf C}({\mathbf T}))\cong  K^*_W({\mathbf
T}) \cong K_*(C^*(G))\cong$$
$$\cong \HE_*(C^*(G))\cong \HE^W_*({\mathbf C}({\mathbf T})\cong
H^W_*({\mathbf C}({\mathbf T})\cong H^*_W({\mathbf
T}) \cong H^*(G) \cong \HP_*(C^*(G))$$ are
as follows.

(a). For any compact Lie group $G$, let $R[G]$ be the representation ring.
Then
the
${\mathbf Z}/(2)$-graded algebra $K^*(G) =
\wedge_{\mathbf
C}(\beta(\rho_1), \beta(\rho_2),\dots,\beta(\rho_n))$, where $\rho_i$ are
the standard
irreducible representations and $\beta :R[G] \to K^*(G)$ is the Bott map.
Hence, from
we have $$K^*(\SU(n+1)) \cong \wedge_{\mathbf C}(\beta(\rho_1),
\dots,
\beta(\rho_n)),$$ $$K^*(\SO(2n+1)) \cong \wedge_{\mathbf C}(\beta(\rho_1),
\dots,
\beta(\rho_n), \varepsilon_{2n+1}).$$

(b). It follows that the ${\mathbf Z}/(2)$-graded complex
cohomology
groups are exterior algebras over ${\mathbf C}$ and in particular
$$H^*(\SU(2n)) \cong \wedge_{\mathbf C}(x_3,x_5,\dots, x_{4n-1}),$$
$$H^*(Sp(n))
\cong \wedge_{\mathbf C}(x_3,x_7,\dots, x_{4n-1}),$$ $$H^*(\SU(2n+1)) \cong
\wedge_{\mathbf C}(x_3,x_5,\dots, x_{4n+1}),$$ $$H^*(\SO(2n+1)) \cong
\wedge_{\mathbf C}(x_3,x_7,\dots, x_{4n-1}).$$

(c). Define a function $\Phi: {\mathbf N} \times {\mathbf N} \times
{\mathbf N} \to
{\mathbf Z}$ by
$$\Phi(n,k,\varepsilon) = \sum_{i=1}^k (-1)^{i-1} \binom{n}{k-i} i^{q-1}.$$
It then
follows that we have Chern character $ch:
K^*(\SU(n+1)) \to H^*(SU(n+1))$ given by
$$ch(\beta(\rho_k)) =
\sum_{i=1}^n \frac{(-1)^i}{i!}\Phi(n+1,k,i+1) x_{2i+1}, \forall k \geq 1,$$
$ch : K^*(\SO(2n+1)) \to
H^*(\SO(2n+1)),$ given by the formula $$ch(\beta(\lambda_k)) =
\sum_{i=1}^n
\frac{(-1)^{i-1}2}{(2i-1)!}\Phi(2n+1,k,2i)x_{4i-1} (\forall k=1,2,\dots,n-1)$$
$$ch(\varepsilon_{2n+1}) = \sum_{i=1}^n
\frac{(-1)^{i-1}}{2^{n-1}(2i-1)!}\sum_{k=1}^n\Phi(2n+1,k,2i)x_{4i-1}.$$
\end{remark}

\section{Chern Characters of Compact Quantum Groups}
{\bf (joint with A. O. Kuku and N. Q. Tho)}

The unitary representations of compact quantum groups
are classified as follows.

\begin{itemize}
\item Every irreducible unitary representation of ${\mathcal
F}_\varepsilon(G)$ on a Hilbert space  is a completion of a unitarizable
highest weight representation. Moreover, two such representations are
equivalent if and only if they have the same weight.
\item The highest weight representations can be described as follows. Let 
$w = s_{i_1}.s_{i_2}. \dots s_{i_k}$ be a reduced decomposition of an
element $w$ of the Weyl group $W$. Then, (i) the Hilbert space tensor
product 
$$\rho_{w,t} = \pi_{s_{i_1}} \otimes \dots \otimes \pi_{s_{i_k}} \otimes
\tau_t$$ is an irreducible *-representation of ${\mathcal
F}_\varepsilon(G)$ which is associated to the Schubert cell ${\mathbf
S}_w$; (ii) up to equivalence, the representation $\rho_{w,t}$ does not
depend on the choice of the reduced decomposition of $w$; (iii) every
irreducible *-representation of ${\mathcal F}_\varepsilon(G)$ is
equivalent to some $\rho_{w,t}$. 
\end{itemize}

Moreover, one can show that 
$$\bigcap_{(w,t)\in W \times {\mathbf T}}\ker \rho_{w,t} = \{ e \},$$ i.e.
the
representation 
$$\bigoplus_{w\in W} \int_{\mathbf T}^\oplus \rho_{w,t} dt$$ is faithful
and 
$$\dim \rho_{w,t} = \left\{ \begin{array}{cl} 1, & \qquad \mbox{ if } w =
e,\\
				             \infty, & \qquad \mbox{ if } 
w \ne e
\end{array} \right. $$

We record now the definition of C*-algebraic compact quantum group. 
Let $G$ be a complex algebraic group with compact real form.
\begin{definition}The {\sl C*-algebraic
compact
quantum group} $C^*_\varepsilon(G)$ is the C*-completion
 of the *-algebra ${\mathcal F}_\varepsilon(G)$ with respect to the
C*-norm $$\Vert f \Vert := \sup_\rho \Vert \rho(f) \Vert \quad (f\in
{\mathcal F}_\varepsilon(G)),$$ where $\rho$ runs through the
*-representations of ${\mathcal F}_\varepsilon(G)$ and the norm on the
right-hand side is the operator norm.  \end{definition}

We now prove the following result about the structure of compact quantum
groups.
\begin{theorem} 
$$C^*_\varepsilon(G) \cong {\mathbf C}({\mathbf T}) \oplus \bigoplus_{e
\ne w\in
W}\int_{\mathbf T}^{\oplus} {\mathcal K}({\mathbf H}_{w,t}) dt,$$
where ${\mathbf C}({\mathbf T})$ is the algebra of complex valued continuous
function on ${\mathbf
T}$ and ${\mathcal K}({\mathbf H})$ is the ideal of compact operators in a separable
Hilbert space ${\mathbf H}$ 
\end{theorem}

Let $A$ be a involutive Banach  algebra. 
We construct a non-commutative character $$ch_{C^*} : K_*(A) \to \HE_*(A)$$ and
later show
that when $A= C^*_\varepsilon(G)$, this Chern character reduces up to isomorphism to
classical Chern character on the normalizers of maximal compact tori. 

Let $A$ be an involutive Banach algebra with unity. 
\begin{theorem}
There exists a Chern character $$ch_{C^*} : K_*(A) \to \HE_*(A).$$ 
\end{theorem}

Our next result computes the Chern characters  for $A=C^*_\varepsilon(G)$ by
reducing it to the classical case.

\begin{theorem}
Let ${\mathbf T}$ be a fixed maximal torus of $G$ with Weyl group $W:=
{\mathcal N}_{\mathbf T}/{\mathbf T}$.
Then the Chern character $$ch_{C^*}: K_*(C^*_\varepsilon(G)) \to
\HE_*(C^*_\varepsilon(G))$$ is an
isomorphism modulo torsion, i.e.
$$\CD
ch_{C^*}: K_* (C^*_\varepsilon(G))\otimes {\mathbf C} @>\cong >>
\HE_*(C^*_\varepsilon(G)),\endCD $$
which can be identified with the classical
Chern character
$$\CD ch: K_*({\mathbf C}({\mathcal N}_{\mathbf T})) @>>> \HE_*({\mathbf
C}({\mathcal N}_{\mathbf T}))\endCD$$
that
is also an isomorphism modulo torsion, i.e.
$$\CD
ch: K_* ({\mathcal N}_{\mathbf T})\otimes {\mathbf C} @>\cong >>
H^*_{DR}({\mathcal N}_{\mathbf T}).
\endCD$$ 
\end{theorem}

\newpage
{\bf References}
\vskip 1cm
\par
1. Do Ngoc Diep, \textit{A survey of Noncommutative geometry
methods for group algebras}, {\bf Journal of Lie Theory},
\textbf{3}(1993), 149-177.

\par
2. Do Ngoc Diep, {\it Noncommutative geometry Methods for Group
C*-Algebras},282pp., to appear in {\bf Pitman Mathematics Series of
Research
Notes}; \\ http://xxx.lanl.gov/math.KT/9807124

\par
3. D. N. Diep and N. V. Thu', \textit{Homotopy invariance of
entire current
cyclic homology}, {\bf Vietnam J. of Math. (Springer)}, \textbf{25}(1997),
No 3, 211-228.

\par 
4. D. N. Diep and N. V. Thu', \textit{Entire homology of
noncommutative de Rham currents}, to appear in {\bf Publication of Center 
for Complex Analysis and Functional Anlysis}, Vietnam
National Univ; ICTP Preprint IC96/214. 

\par
5. D. N. Diep, A. O. Kuku and N. Q. Tho,
\textit{Noncommutative Chern
characters of compact Lie group C*-algebras},{\bf K-Theory} (to appear);\\
http://xxx.lanl.gov/math.KT/9807101 .

\par
6. D. N. Diep, A. O. Kuku and N. Q. Tho,
\textit{Noncommutative Chern
characters of compact quantum groups
},19pp. \\ http://xxx.lanl.gov/math.KT/9807099.

\vskip 1cm
{\bf Address:} Institute of Mathematics, National Center for Science and
Technology, P. O. Box 631, Bo Ho, VN-10.000, Hanoi, Vietnam

{\bf Email:} {\tt dndiep@@ioit.ncst.ac.vn; \\ dndiep@@member.ams.org
(forward mail)}

\end{document}